\documentclass[12pt]{article}

\usepackage{amsthm}
\usepackage{amsmath}
\usepackage{amsfonts}
\usepackage{amssymb}
\usepackage{amscd}
\usepackage{mathrsfs}
\usepackage[all]{xypic}
\usepackage{url}

\newtheorem{thm}{Theorem}[section]

\newtheorem{prop}[thm]{Proposition}

\theoremstyle{definition}

\theoremstyle{definition}

\newcommand{\sinpioverthree}{0.866025404}
\newcommand{\C}{\mathbb{C}}

\newcommand{\R}{\mathbb{R}}

\newcommand{\indep}{\! \perp \!\!\! \perp \!}

\newcommand{\szero}{0}
\newcommand{\splus}{+ \! \!\!\! +}

\newcommand{\cia}[3]{\ensuremath{ [ #1 \! \indep \! #2 | #3 ] }}
\newcommand{\ciaa}[3]{\ensuremath{ [\![ #1 \! \indep \! #2 | #3 ]\!] }}

\begin{document}
\bibliographystyle{plain}
\title{Three Counterexamples on Semigraphoids}
\author{Raymond Hemmecke, Jason Morton, Anne Shiu, \\Bernd Sturmfels, and Oliver Wienand}
\date{}

\maketitle
\abstract{Semigraphoids are combinatorial structures that arise in
statistical learning theory. They are equivalent to convex rank
tests and to polyhedral fans that coarsen the reflection arrangement
of the symmetric group $S_n$. In this paper we resolve two problems
on semigraphoids posed in Studen\'y's book
\cite{Studeny2005Probabilistic}, and we answer a related question 
 by Postnikov,  Reiner, and Williams on
generalized permutohedra \cite{PRW}.
We also study the semigroup and the toric ideal associated with
semigraphoids.}

\section{Introduction}

A  {\em conditional independence (CI) statement} on a finite set of
random variables, indexed by $[n] = \{1,2,\ldots,n\}$,  is a formal
symbol $\, [ i \indep j \, |\, K]\,$ where $K \subset [n]$ and $i,j
\in [n] \backslash K$. The symbol $\, [ i \indep j \, |\, K]\,$
represents the statement that the random variables $i$ and $j$ are
conditionally independent given the joint random variable $K$. For
any joint probability distribution on the $n$ random variables, the
set $\mathcal{M}$ of all  CI statements that are valid for the given
distribution satisfies the following axiom: \break

\noindent ({\bf SG}) If $  [i \indep j \, | K \cup {\ell}]$ and $ [i
\indep \ell \, | K]$ are in $ \mathcal{M}$ then so are $ [i \indep j
\,| K]\,$ and $[i \indep \ell \,|  K \cup j] $. \break

\noindent A {\em semigraphoid} is any set $\mathcal{M}$ of CI
statements which satisfies the axiom ({\bf SG}). Studen\'y's book
\cite{Studeny2005Probabilistic} gives an introduction to
semigraphoids and their role in statistical learning theory. For
further details and references see also Mat\' u\v s  \cite{Mat1,
Mat3}. In this paper we construct examples which answer two problems stated by Studen\'y:

\vskip .07cm

\noindent ({\bf Q1}) {\em Is it true that every coatom of the
lattice of (disjoint) semigraphoids over $[n]$ is a structural
independence model over $[n]$?} \ \ \  \phantom{dadada} \,\qquad
\cite[Question 4, page 194]{Studeny2005Probabilistic} \ \ \ \break
\noindent ({\bf Q2}) {\em Is every structural imset over $[n]$
already a combinatorial imset over $[n]$?} \ \ \
\phantom{dadadadadadadada dadadadadadada dadadadadad} \qquad
\cite[Question 7, page 207]{Studeny2005Probabilistic}

\smallskip

Our approach is based on the geometric characterization of
semigraphoids which was developed in  \cite{GRT}. Let $\Pi_{n-1}$
denote the $(n \!- \! 1)$-dimensional {\em permutohedron}
\cite{Mat2, Zie}, and let $C_n = [0,1]^n$ denote the standard
$n$-dimensional cube. The vertices of $\Pi_{n-1}$ are in bijection
with the elements of the symmetric group $S_n$, and with the
monotone edge paths from $(0,0,\ldots,0)$ to $(1,1,\ldots,1)$ on the
cube $C_n$. The $2$-dimensional faces of $C_n$ are in bijection with
the CI statements on $[n]$. Namely, $\,[ i \indep j \, |\, K ] = [ j
\indep i \, |\, K ]\,$ represents the $2$-face of $C_n$ with
 $x_k = 1$ for $k \in K$ and $x_l = 0$ for $l \in [n]
\backslash (K \cup \{i,j\}) $. The number of these $2$-cubes equals
$\,\gamma_n := \, \binom{n}{2} 2^{n-2} $.
There is a natural surjection from the edges of $\Pi_{n-1}$ onto the
$2$-faces of $C_n$. Namely, an edge of $\Pi_{n-1}$ corresponds to a
pair of adjacent monotone edge paths on $C_n$. These adjacent paths
differ only along a $2$-cube $[ i \indep j \, |\, K ]$. In this manner,
we identify any set $\mathcal{M}$ of CI statements on $[n]$ with
a set of $2$-cubes on the boundary of $C_n$.
We also identify $\mathcal{M}$ with a set of edges of the permutohedron
$\Pi_{n-1}$, bearing in mind that opposite edges of a square have the same CI statement as their label.

Each $2$-face of the permutohedron
 $\Pi_{n-1}$ is either a square or a hexagon. By \cite{GRT},
the semigraphoid axiom is equivalent to the following geometric
condition on $\Pi_{n-1}$:

\vskip .04cm \noindent (${\bf SG}'$) {\em If two adjacent edges of a
hexagon are in $\mathcal{M}$ then so are their two opposites.}
\vskip .04cm
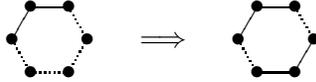
\[
\begin{xy}<5mm,0cm>:
(-.5,\sinpioverthree)  ="123"   *{\bullet};
(.5 ,\sinpioverthree)  ="132"   *{\bullet};
(1,0)                  ="312"   *{\bullet};
(.5,-\sinpioverthree)  ="321"   *{\bullet};
(-.5,-\sinpioverthree) ="231"   *{\bullet};
(-1,0)                 ="213"   *{\bullet};
   "123";"132" **@{-};
   "132";"312" **@{.};
   "312";"321" **@{.};
   "321";"231" **@{.};
   "231";"213" **@{.};
   "213";"123" **@{-};
\end{xy} \quad
\implies \quad
\begin{xy}<5mm,0cm>:
(-.5,\sinpioverthree)  ="123"   *{\bullet};
(.5 ,\sinpioverthree)  ="132"   *{\bullet};
(1,0)                  ="312"   *{\bullet};
(.5,-\sinpioverthree)  ="321"   *{\bullet};
(-.5,-\sinpioverthree) ="231"   *{\bullet};
(-1,0)                 ="213"   *{\bullet};
   "123";"132" **@{-};
   "132";"312" **@{.};
   "312";"321" **@{-};
   "321";"231" **@{-};
   "231";"213" **@{.};
   "213";"123" **@{-};
\end{xy}
\]
The normal fan of the permutohedron $\Pi_{n-1}$ is
the reflection arrangement of $S_n$. 
Theorem 3 in \cite{GRT} identifies semigraphoids with fans that
coarsen this arrangement.
Such fans are called {\em
convex rank tests}. Namely, $\mathcal{M}$ specifies the set of edges
of $\Pi_{n-1}$ whose dual walls in the normal fan are not
present in
 the convex rank test.

A basic question about any semigraphoid $\mathcal{M}$ is whether its corresponding
convex rank test is {\em submodular}, in other words, whether it is the normal fan of a
convex polytope. That polytope would then be a Minkowski summand of
$\Pi_{n-1}$. These polytopes are known as {\em generalized permutohedra}
and they were studied in
\cite{Postnikov2005, PRW}. 

Studen\'y's first question has the
following geometric translations:

\vskip .03cm \noindent  ({\bf Q1}) {\em Is every coarsest convex
rank test submodular?}  \hfill \break \noindent  ({\bf Q1}) {\em Is
every fan which maximally coarsens the $S_n$-arrangement the normal
fan of a generalized permutohedron?} \vskip .03cm

\noindent
In the first version of \cite{PRW},
Postnikov, Reiner and Williams asked a 
similar question:

\noindent  ({\bf Q3}) {\em Is every simplicial fan which coarsens
the $S_n$-arrangement the normal fan of a simple generalized
permutohedron?} \vskip .03cm

This paper answers all three questions. In Section 2 we derive and explain our
counterexample for Question ({\bf Q3}). That example is discussed in
\cite[Example 3.8]{PRW}. 
 By Studen\'y's classification of the $26424$ semigraphoids for $n=4$, it had been known
that the answers to Questions ({\bf Q1}) and ({\bf Q2}) are affirmative for $n
\leq 4$. In Sections 3 and 4 we construct counterexamples for 
({\bf Q1}) and ({\bf Q2}) with $n=5$.

Question ({\bf Q2}) has the following reformulation in the setting
of toric algebra \cite[\S 7]{MS}. We represent the semigraphoid
axiom as an equation among formal symbols:
$$ ({\bf SG}'') \qquad \qquad
  [i \indep j \, | K \cup {\ell}] \,+ \, [i \indep \ell \, | K] \quad = \quad
     [i \indep j \,| K]\, + \, [i \indep \ell \,|  K \cup j]  \qquad \qquad $$
for all $i,j,l,K$.  These relations span the kernel of the linear map
\begin{equation}
\label{Amap} \,\mathcal{A}\, :\, \mathbb{Z}^{\gamma_n} \,\rightarrow
\,\mathbb{Z}^{2^n}\,,\,\,\, [i \indep j \,| \, K] \,\,\, \mapsto
\,\,\, 
 e_{iK} + e_{jK} - e_{K} - e_{ijK} .
\end{equation}
A semigraphoid is a solution to the equations $({\bf SG}'')$  in the
semiring $\{\szero, \splus \}$, representing
``zero'' and ``positive''.
A semigraphoid is {\em submodular} if it is the set of 
zero coordinates of a solution to $({\bf SG}'')$
 in the  non-negative real numbers.
These definitions furnish us with an algebraic
representation of a semigraphoid $\mathcal{M}$ and a systematic
method for testing submodularity of $\mathcal{M}$ by linear
programming. Studen\'y's question ({\bf Q2}) concerns the 
$\mathbb{N}$-linear span of
the columns of the matrix~$\mathcal{A}$:

\vskip .04cm \noindent ({\bf Q2}) {\em Is the semigroup
$\mathcal{A}( \mathbb{N}^{\gamma_n})$ normal, i.e., does it coincide
with $\, \mathcal{A}( \mathbb{R}_{\geq 0}^{\gamma_n} )\,\cap \,
\mathbb{Z}^{2^n}$~?} \vskip .04cm

In Section 5 we study the toric ideal \cite{ATY} of $\,\mathcal{A}\,$ in a
polynomial ring in $\gamma_n$ unknowns, and we examine how it
differs from the subideal generated by the binomials
$$ ({\bf SG}''')  \qquad \qquad \quad
  [i \indep j \, | K \cup {\ell}] \,\cdot \, [i \indep \ell \, | K] \,\,\, - \,\,\,
     [i \indep j \,| K]\, \cdot \, [i \indep \ell \,|  K \cup j] . \qquad \qquad \quad .$$
  Proposition  \ref{primedec} describes the primary decomposition
 of this binomial ideal for $n = 4$. We also discuss the problem of 
 deriving the full Markov basis from (${\bf SG}'''$).

\section{A non-submodular simplicial semigraphoid}

Let $n=4$ and consider the $4$-dimensional cube $C_4$ and the
$3$-dimensional permutohedron $\Pi_3$. Each hexagon on  $\Pi_3$
corresponds to one of the eight facets of $C_4$. Each facet
specifies three semigraphoid axioms, written additively as in $({\bf
SG}'')$:
$$ \begin{matrix}
&& \ciaa{1}{2}{\emptyset} + \cia{2}{3}{1} \, &=& \,\cia{2}{3}{\emptyset} + \cia{1}{2}{3} & \Longleftarrow \\
(*,*,*,0) && \cia{1}{3}{\emptyset} + \cia{1}{2}{3} \, &=& \,\ciaa{1}{2}{\emptyset} + \cia{1}{3}{2} \\
&& \cia{1}{3}{\emptyset} + \cia{2}{3}{1} \, &=&
\,\cia{2}{3}{\emptyset} + \cia{1}{3}{2}
\end{matrix} $$  \smallskip  $$ \begin{matrix}
&& \ciaa{1}{2}{\emptyset} + \cia{2}{4}{1} \, &=& \,\cia{2}{4}{\emptyset} + \cia{1}{2}{4} & \quad \;\;\; \\
(*,*,0,*) && \ciaa{1}{2}{\emptyset} + \cia{1}{4}{2} \, &=& \,\cia{1}{4}{\emptyset} + \cia{1}{2}{4} \\
&& \cia{1}{4}{\emptyset} + \cia{2}{4}{1} \, &=&
\,\cia{2}{4}{\emptyset} + \cia{1}{4}{2}
\end{matrix} $$  \smallskip   $$ \begin{matrix}
&& \cia{1}{3}{\emptyset} + \cia{1}{4}{3} \, &=& \,\cia{1}{4}{\emptyset} + \cia{1}{3}{4}  \\
(*,0,*,*) && \ciaa{3}{4}{\emptyset} + \cia{1}{3}{4} \, &=& \,\cia{1}{3}{\emptyset} + \cia{3}{4}{1} \\
&& \ciaa{3}{4}{\emptyset} + \cia{1}{4}{3} \, &=&
\,\cia{1}{4}{\emptyset} + \cia{3}{4}{1} & \Longleftarrow \,
\end{matrix} $$  \smallskip  $$ \begin{matrix}
&&\; \cia{2}{3}{\emptyset} + \cia{3}{4}{2} \, &=& \,\ciaa{3}{4}{\emptyset} + \cia{2}{3}{4} & \quad \;\;\; \\
(0,*,*,*) && \cia{2}{4}{\emptyset} + \cia{2}{3}{4} \, &=& \,\cia{2}{3}{\emptyset} + \cia{2}{4}{3} \\
&& \ciaa{3}{4}{\emptyset} + \cia{2}{4}{3} \, &=&
\,\cia{2}{4}{\emptyset} + \cia{3}{4}{2}
\end{matrix} $$ \smallskip  $$ \begin{matrix}
&& \cia{3}{4}{1} + \ciaa{2}{3}{14} \, &=& \,\cia{2}{3}{1} +
\cia{3}{4}{12}
& \Longleftarrow \, \\
(*,*,*,1) && \cia{2}{4}{1} + \ciaa{2}{3}{14} \, &=& \,\cia{2}{3}{1} + \cia{2}{4}{13} \\
&& \cia{2}{4}{1} + \cia{3}{4}{12} \, &=& \,\cia{3}{4}{1} + \cia{2}{4}{13} \\
\end{matrix} $$  \smallskip   $$ \begin{matrix}
&& [1 \!\!\perp\!\!\!\perp\!\!3|2] + [3
\!\!\perp\!\!\!\perp\!\!4|12]\,
 &=& \,[3 \!\!\perp\!\!\!\perp\!\!4|2] + [1 \!\!\perp\!\!\!\perp\!\!3|24]  & \quad \;\;\; \\
(*,*,1,*) && [1 \!\!\perp\!\!\!\perp\!\!3|2] + [\![1
\!\!\perp\!\!\!\perp\!\!4|23]\!]\,
 &=& \,[1 \!\!\perp\!\!\!\perp\!\!4|2] + [1 \!\!\perp\!\!\!\perp\!\!3|24] \\
&& [3 \!\!\perp\!\!\!\perp\!\!4|2] + [\![1
\!\!\perp\!\!\!\perp\!\!4|23]\!]\,
 &=& \,[1 \!\!\perp\!\!\!\perp\!\!4|2] + [3 \!\!\perp\!\!\!\perp\!\!4|12]
\end{matrix} $$  \smallskip   $$ \begin{matrix}
&& [1 \!\!\perp\!\!\!\perp\!\!2|3] + [\![1
\!\!\perp\!\!\!\perp\!\!4|23]\!]\,
 &=& \,[1 \!\!\perp\!\!\!\perp\!\!4|3] + [1 \!\!\perp\!\!\!\perp\!\!2|34]
& \Longleftarrow \\
(*,1,*,*) && [1 \!\!\perp\!\!\!\perp\!\!4|3] + [2
\!\!\perp\!\!\!\perp\!\!4|13]\, &=& \,[2
\!\!\perp\!\!\!\perp\!\!4|3] + [\![1
\!\!\perp\!\!\!\perp\!\!4|23]\!]
 \\
&& [1 \!\!\perp\!\!\!\perp\!\!2|3] + [2
\!\!\perp\!\!\!\perp\!\!4|13]\,
 &=& \,[2 \!\!\perp\!\!\!\perp\!\!4|3] + [1 \!\!\perp\!\!\!\perp\!\!2|34]
\end{matrix} $$ \smallskip  $$ \begin{matrix}
&& [1 \!\!\perp\!\!\!\perp\!\!3|4] + [\![2
\!\!\perp\!\!\!\perp\!\!3|14]\!]\,
 &=& \,[2 \!\!\perp\!\!\!\perp\!\!3|4] + [1 \!\!\perp\!\!\!\perp\!\!3|24] & \quad \;\;\; \\
(1,*,*,*) && [1 \!\!\perp\!\!\!\perp\!\!2|4] + [1
\!\!\perp\!\!\!\perp\!\!3|24]\,
 &=& \,[1 \!\!\perp\!\!\!\perp\!\!3|4] + [1 \!\!\perp\!\!\!\perp\!\!2|34] \\
&& [1 \!\!\perp\!\!\!\perp\!\!2|4] + [\![2
\!\!\perp\!\!\!\perp\!\!3|14]\!]\,
 &=& \,[2 \!\!\perp\!\!\!\perp\!\!3|4] + [1 \!\!\perp\!\!\!\perp\!\!2|34].\\
\end{matrix}
$$
This is a system of $ 24 $ equations in $\,\gamma_4 =
24$ formal symbols $[i \indep j \,|\,K ]$.

A semigraphoid is a solution to these equations over the semiring
$\,\{ \szero, \splus\}$. More precisely, given such a solution vector in
$\{\szero,\splus\}^{24}$, the semigraphoid $\mathcal{M}$ consists of all
coordinates $\,[i \indep j \,|\,K ]\,$ that have the value
$\szero$. There are $26424$ such semigraphoids. They form a sublattice of
the Boolean lattice $\{\szero,\splus\}^{24}$, with $\splus < \szero $.  Question
({\bf Q1}) concerns the coatoms of this lattice. But let us first
resolve Question ({\bf Q3}).

We consider the following collection of CI statements:
 \begin{equation}
\mathcal{M} \,\,\, = \,\,\, \bigl\{\,
 [\![ 2 \perp \!\!\! \perp 3\, |\, 14]\!] ,\,
 [\![1 \perp \!\!\! \perp 4 \,|\,  23 ]\!], \,
 [\![1 \perp \!\!\! \perp 2 \,|\, \emptyset ]\!],\,
[\![ 3 \perp \!\!\! \perp 4 \,| \, \emptyset ]\!] \,\bigr\}.
 \end{equation}
These four symbols are highlighted
 in the $24$ equations above by the use of
double brackets $\, [\![ \,\, \cdots \,\, ]\!]$. Each equation
(individually) can be solved 
among the positive reals after
these four symbols have been set to zero, or equivalently 
 they can be solved as a system over $\{ \szero, \splus\}$. This shows that
$\mathcal{M}$ is a semigraphoid.

\vskip -0.1cm
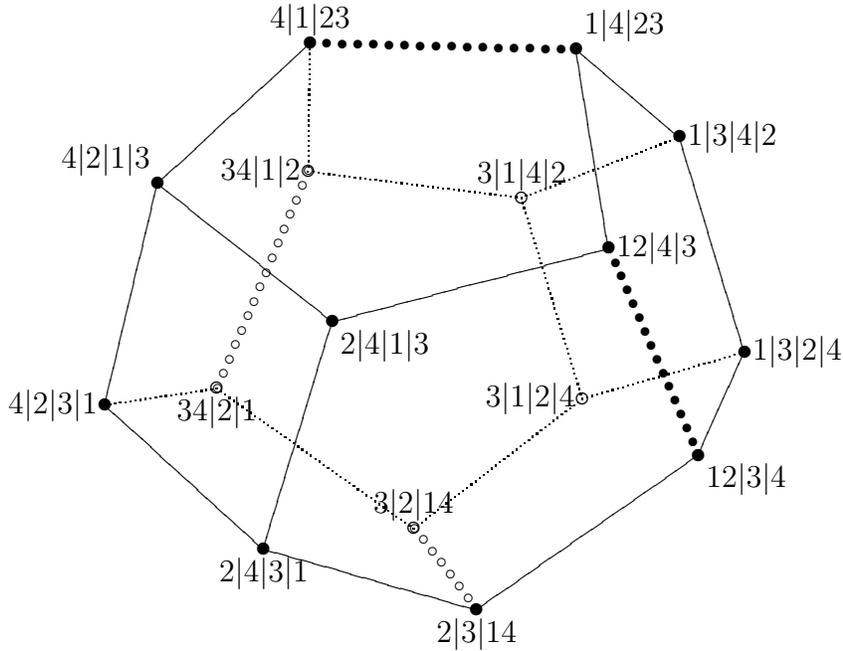
\begin{figure}[htb]
\[
 \begin{xy}<27mm,0cm>:
(2.12,.55)  ="3214-3241"  *+!D{3|2|14} *{\circ}; (2.43,.15)
="2314-2341"  *+!U{2|3|14} *{\bullet};
(3.52,.91)  ="1234-2134"  *+!UL{12|3|4} *{\bullet}; (3.08,1.93)
="1243-2143"  *+!L{12|4|3} *{\bullet};
(1.60,2.31)  ="3412-4312"  *+!R{34|1|2} *{\circ}; (1.15,1.24)
="3421-4321"  *+!U{34|2|1} *{\circ};
(1.61,2.94)  ="4123-4132"  *+!D{4|1|23} *{\bullet}; (2.92,2.91)
="1423-1432"  *+!DL{1|4|23} *{\bullet};
(2.95,1.19)  ="3124"  *+!R{3|1|2|4} *{\circ}; (1.38,.45)  ="2431"
*+!U{2|4|3|1} *{\bullet}; (3.75,1.42)  ="1324"  *+!L{1|3|2|4}
*{\bullet}; (2.65,2.18)  ="3142"  *+!D{3|1|4|2} *{\circ}; (1.72,1.57
)  ="2413"  *+!UL{2|4|1|3} *{\bullet}; (.60,1.16)  ="4231"
*+!R{4|2|3|1} *{\bullet}; (3.43,2.48)  ="1342"  *+!L{1|3|4|2}
*{\bullet}; (.86,2.25)  ="4213"  *+!DR{4|2|1|3} *{\bullet};
"3214-3241";"2314-2341" **@{o};
"1234-2134";"1243-2143" **@{*};
"4123-4132";"1423-1432" **@{*};
"3412-4312";"3421-4321" **@{o};
"4213";"2413" **@{-}; 
"4231";"2431" **@{-}; 
"4213";"4231" **@{-}; 
"2413";"2431" **@{-}; 
"1342";"1324" **@{-}; 
"3142";"3124" **@{.}; 
"1342";"3142" **@{.}; 
"1324";"3124" **@{.}; 
"2314-2341";"1234-2134" **@{-}; 
"3124";"3214-3241" **@{.}; 
"3214-3241";"3421-4321" **@{.}; 
"3412-4312";"3142" **@{.}; 
"1324";"1234-2134" **@{-}; 
"1423-1432";"1342" **@{-}; 
"4123-4132";"3412-4312" **@{.}; 
"1423-1432";"1243-2143" **@{-}; 
"2314-2341";"2431" **@{-}; 
"3421-4321";"4231" **@{.}; 
"2413";"1243-2143" **@{-}; 
"4123-4132";"4213" **@{-}; 
\end{xy}
\]
\vskip -0.2cm
\caption{A simple 3-dimensional polytope with $16$ vertices and $10$
facets} \label{nonpolytopalpolytope}
\end{figure}

The semigraphoid $\mathcal{M}$ is represented geometrically
 by the three-dimensional polytope in Figure 1.
This polytope is {\em simple}, i.e., each of the $16$ vertices is
adjacent to three other vertices. The eight vertices whose
labels include three bars (such as $4|2|1|3$) correspond to unique
permutations in $S_4\,$ (namely the permutation $4213$), while the
eight vertices whose labels have two bars (such as $4|1|23$)
correspond to pairs of permutations in $S_4\,$ (namely $4123$ and
$4132$). This partition of $S_4$ into eight singletons and eight
pairs is the convex rank test of $\mathcal{M}$. The
normal fan of the polytope in Figure 1 is a simplicial fan which
is combinatorially (but not geometrically) isomorphic to a fan that
coarsens the hyperplane arrangement of $S_4$.

\begin{prop} \label{mmprop}
The simplicial semigraphoid $\mathcal{M}$ is not submodular.
\end{prop}

\begin{proof}
Suppose that $\mathcal{M}$ were submodular. Then the above equations
have a solution in $(\mathbb{R}_{\geq 0})^{24}$ whose coordinates in
$\mathcal{M}$ are zero and whose other $20$ coordinates are
positive. The four equations marked by an ``$\Longleftarrow$'' give
the following four equations:
$$
\begin{matrix}

[2 \!\!\perp\!\!\!\perp\!\!3|1]\, &=& \,[2
\!\!\perp\!\!\!\perp\!\!3| \emptyset ]
+ [1 \!\!\perp\!\!\!\perp\!\!2|3]  \\
 [1 \!\!\perp\!\!\!\perp\!\!4|3]\,
 &=& \,[1 \!\!\perp\!\!\!\perp\!\!4|  \emptyset ] + [3 \!\!\perp\!\!\!\perp\!\!4|1]
\\
\,[3 \!\!\perp\!\!\!\perp\!\!4|1]
 &=&
[2 \!\!\perp\!\!\!\perp\!\!3|1] + [3 \!\!\perp\!\!\!\perp\!\!4|12]\,
\\
[1 \!\!\perp\!\!\!\perp\!\!2|3]
 &=& \,[1 \!\!\perp\!\!\!\perp\!\!4|3] + [1 \!\!\perp\!\!\!\perp\!\!2|34].
\end{matrix}
$$
Adding the left hand sides and the right hand sides of the four equations yields
$$
[2  \!\!\perp\!\!\!\perp\!\!3|  \emptyset  ]\, + \, [1
\!\!\perp\!\!\!\perp\!\!4|  \emptyset  ] \,+\,
 [3 \!\!\perp\!\!\!\perp\!\!4 |12] \,+\,
 [1 \!\!\perp\!\!\!\perp\!\!2|34] \quad = \quad 0.
$$
This contradicts the assumption that these four values are strictly
positive.  \end{proof}

The set of all non-negative solutions to the $24$ equations is an
$11$-dimensional cone in $(\mathbb{R}_{\geq 0})^{24}$.  This cone is isomorphic to the
$16$-dimensional cone of submodular functions on $2^{[4]}$, modulo its $5$-dimensional lineality space. Its $22108$ faces are in
bijection with the submodular semigraphoids, or, equivalently, with
the generalized permutohedra for $n=4$. In addition to these, there
are $   4316$ semigraphoids that are not submodular. Each of the latter
can be represented by a polytope of dimension $\leq 3$  as in Figure
1. These polytopes have the combinatorial properties  of
generalized permutohedra, but they cannot be realized as
Minkowski summands of $\Pi_3$. For example, see \cite[Figure
5]{Hirai} for a polytope that depicts Studen\'{y}'s example of a
semigraphoid that is not submodular (see \cite{GRT} and
\cite[Section 2.2.4]{Studeny2005Probabilistic}). 

We now give a classification of non-submodular semigraphoids for $n=4$
and $|\mathcal{M}$ small. All simplicial examples
 are coarsenings (up to relabeling) of the
 particular semigraphoid $\mathcal{M}$
 in Proposition \ref{mmprop}.
   The following table lists the number of semigraphoids classified by number of CI statements, their type, and whether they are simplicial.  Here,
the {\em type} of a semigraphoid is the triple $(m_0,m_1,m_2)$ where $m_t$ is the number of CI statements $\, [ i \indep j \, |\, K]\,$ in $\mathcal{M}$ such that $|K| = m_t$.
\bigskip

\begin{center}
  \begin{tabular}{ | l || c | c | c | c | }
    \hline
    $| \mathcal{M} |$ & type & non-simplicial &  simplicial & total\\ \hline
 3 & ( 0 , 3 , 0 ) & 8 & 0 & 8  \\ \hline
4 & ( 0 , 4 , 0 ) & 78 & 0 & 78  \\ \hline
4 & ( 1 , 2 , 1 ) & 30 & 0 & 30 \\ \hline
4 & ( 2 , 0 , 2 ) & 0 & 6 & 6 \\ \hline
5 & ( 0 , 5 , 0 ) & 300 & 0 & 300 \\ \hline
5 & ( 1 , 2 , 2 ) & 30 & 0 & 30 \\ \hline
5 & ( 1 , 3 , 1 ) & 84 & 0 & 84 \\ \hline
5 & ( 2 , 0 , 3 ) & 12 & 12 & 24 \\ \hline
5 & ( 2 , 2 , 1 ) & 30 & 0 & 30 \\ \hline
5 & ( 3 , 0 , 2 ) & 24 & 0 & 24 \\ \hline
6 & ( 0 , 6 , 0 ) & 604 & 0 & 604 \\ \hline
6 & ( 1 , 3 , 2 ) & 84 & 0 & 84 \\ \hline
6 & ( 1 , 4 , 1 ) & 78 & 0 & 78 \\ \hline
6 & ( 2 , 0 , 4 ) & 30 & 3 & 33 \\ \hline
 \end{tabular} \end{center} \begin{center}   \begin{tabular}{ | l || c | c | c | c | }
 $| \mathcal{M} |$ & type & non-simplicial &  simplicial & total\\ \hline
6 & ( 2 , 2 , 2 ) & 30 & 0 & 30 \\ \hline
6 & ( 2 , 3 , 1 ) & 84 & 0 & 84 \\ \hline
6 & ( 3 , 0 , 3 ) & 74 & 12 & 96 \\ \hline
6 & ( 4 , 0 , 2 ) & 30 & 3 & 33 \\ \hline
7 & ( 0 , 7 , 0 ) & 684 & 0 & 684 \\ \hline
7 & ( 1 , 4 , 2 ) & 78 & 0 & 78 \\ \hline
7 & ( 1 , 5 , 1 ) & 24 & 0 & 24 \\ \hline
7 & ( 2 , 0 , 5 ) & 18 & 0 & 18 \\ \hline
7 & ( 2 , 3 , 2 ) & 84 & 0 & 84 \\ \hline
7 & ( 2 , 4 , 1 ) & 78 & 0 & 78 \\ \hline
7 & ( 3 , 0 , 4 ) & 132 & 0 & 132 \\ \hline
7 & ( 4 , 0 , 3 ) & 132 & 0 & 132 \\ \hline
7 & ( 5 , 0 , 2 ) & 18 & 0 & 18 \\ \hline
8 & ( 0 , 8 , 0 ) & 450 & 0 & 450 \\ \hline
8 & ( 1 , 5 , 2 ) & 24 & 0 & 24 \\ \hline
8 & ( 2 , 0 , 6 ) & 3 & 0 & 3 \\ \hline
8 & ( 2 , 4 , 2 ) & 48 & 0 & 48 \\ \hline
8 & ( 2 , 5 , 1 ) & 24 & 0 & 24 \\ \hline
8 & ( 3 , 0 , 5 ) & 72 & 0 & 72 \\ \hline
8 & ( 4 , 0 , 4 ) & 174 & 0 & 174 \\ \hline
8 & ( 5 , 0 , 3 ) & 72 & 0 & 72 \\ \hline
8 & ( 6 , 0 , 2 ) & 3 & 0 & 3 \\ \hline
9 & ( 0 , 9 , 0 ) & 212 & 0 & 212 \\ \hline
9 & ( 3 , 0 , 6 ) & 12 & 0 & 12 \\ \hline
9 & ( 4 , 0 , 5 ) & 84 & 0 & 84 \\ \hline
9 & ( 5 , 0 , 4 ) & 84 & 0 & 84 \\ \hline
9 & ( 6 , 0 , 3 ) & 12 & 0 & 12 \\ \hline
10 & ( 0 , 10 , 0 ) & 60 & 0 & 60 \\ \hline
10 & ( 4 , 0 , 6 ) & 15 & 0 & 15 \\ \hline
10 & ( 5 , 0 , 5 ) & 24 & 0 & 24 \\ \hline
10 & ( 6 , 0 , 4 ) & 15 & 0 & 15 \\ \hline
11 & ( 0 , 11 , 0 ) & 12 & 0 & 12 \\ \hline
11 & ( 5 , 0 , 6 ) & 6 & 0 & 6 \\ \hline
11 & ( 6 , 0 , 5 ) & 6 & 0 & 6 \\
  \end{tabular}
\end{center}

\section{A non-submodular coarsest semigraphoid}  \label{sec:coarsest}

We now consider the case $n=5$. There are $\gamma_5 = 80$ CI
statements, one for each two-dimensional face of the  $5$-cube
$C_5$. There are $120$ semigraphoid axioms $({\bf SG}'')$, three
for each of the $40$ three-dimensional faces of $C_5$,
listed as additive equations in the Appendix. The semigraphoids
are the solutions of these equations over $\{\szero,\splus\}^{80}$. 
These solutions include
the all-zero vector ${\bf 0} $ which represents the
semigraphoid that consists of all $80$ CI statements, and which is the
maximal element in the lattice of semigraphoids. A semigraphoid is
said to be {\em coarsest} if it is maximal among non-${\bf 0}$ semigraphoids. Geometrically, such a semigraphoid corresponds to a fan which
coarsens the $S_5$-arrangement but cannot be coarsened to a
non-trivial fan.

We now present the counterexample which answers question ({\bf Q1}).
Our constructions make use of the identification of semigraphoids
with convex rank tests that was derived in \cite{GRT}. Let $\Gamma$
denote the partition of the symmetric group $S_5$ into fourteen
classes as follows. There are eight classes containing $12$
permutations each:
$$
\begin{matrix}
    15|234  &&  234|15 &&    123|45 &&   235|14 \\
    124|35  &&   245|13 &&   134|25 &&  345|12.
    \end{matrix}
$$
And there are six classes containing four permutations each:
$$
\begin{matrix}
   12|5|34   &&  25|1|34  &&  13|5|24  \\
   35|1|24   &&  14|5|23  &&   45|1|23.
   \end{matrix}
$$
Here $\,15|234\,$ denotes the class of all permutations $\,ijklm \,$
with $\{i,j\} = \{1,5\}$ and $\{k,l,m\} = \{2,3,4\}$. Similarly,
$\,45|1|23\,$ denotes the class of all permutations $\,ijklm\,$ with
$\{i,j\} = \{4,5\}$, $k=1$, and $\{l,m\} = \{2,3\}$. Clearly,
$\Gamma$ is a pre-convex rank test, as each of the $14$ classes
is the set of all linear extensions of a poset on $[5] =
\{1,2,3,4,5\}$. Note that the stabilizer of the pre-convex rank test
$\Gamma$ in $S_5$ has order $12$, because $\Gamma$ is fixed under
permutations of $\{ 1,5\}$ and permutations of $\{2,3,4\}$. The $14$
classes of $\Gamma$ are represented by the $14$ vertices of the
polytope in Figure 2.

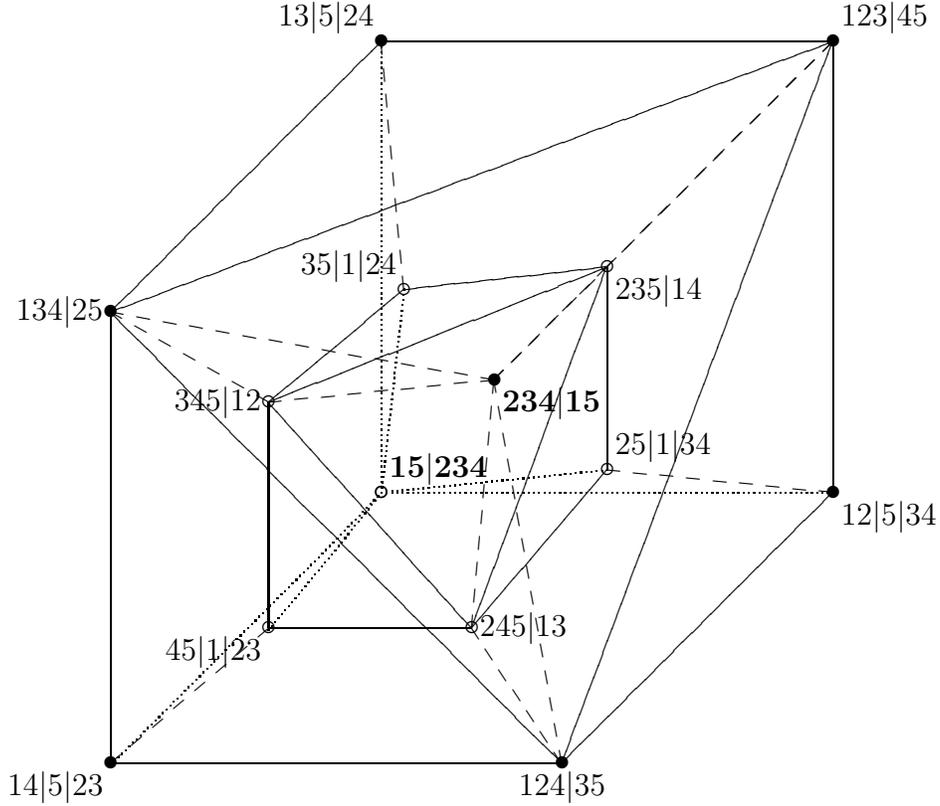
\begin{figure}[htb]
\[
\begin{xy}<30mm,0cm>:
(1.7,1.7) ="Extra" *+!UL{{\bf 234|15}} *{\bullet};
(0,0) ="000" *+!UR{14|5|23} *{\bullet}; (2,0) ="001" *+!U{124|35}
*{\bullet}; (0,2) ="100" *+!R{134|25} *{\bullet};
(1.2,1.2) ="010" *+!LD{{\bf 15|234}} *{\circ}; 
(1.2,3.2) ="110" *+!RD{13|5|24} *{\bullet}; (3.2,1.2) ="011"
*+!UL{12|5|34} *{\bullet}; (3.2,3.2) ="111" *+!LD{123|45}
*{\bullet};
"100";"001" **@{-}; "001";"111" **@{-}; "100";"111" **@{-};
"000";"001" **@{-}; "000";"100" **@{-}; "000";"010" **@{.};
"010";"011" **@{.}; "010";"110" **@{.}; "011";"111" **@{-};
"110";"111" **@{-}; "100";"110" **@{-}; "001";"011" **@{-};
%
(.7,.6) ="i000" *+!UR{45|1|23} *{\circ}; (1.6,.6) ="i001"
*+!L{245|13} *{\circ}; (0.7,1.6) ="i100" *+!R{345|12} *{\circ};
(1.2,1.2) ="i010" *+!RD{} *{\circ}; 
(1.3,2.1) ="i110" *+!RD{35|1|24} *{\circ}; (2.2,1.3) ="i011"
*+!DL{25|1|34} *{\circ}; (2.2,2.2) ="i111" *+!UL{235|14} *{\circ};
"i100";"i001" **@{-}; "i001";"i111" **@{-}; "i100";"i111" **@{-};
"i000";"i001" **@{-}; "i000";"i100" **@{-}; "i000";"i010" **@{.};
"i010";"i011" **@{.}; "i010";"i110" **@{.}; "i011";"i111" **@{-};
"i110";"i111" **@{-}; "i100";"i110" **@{-}; "i001";"i011" **@{-};
%
"Extra";"111" **@{--}; "Extra";"i111" **@{--}; "Extra";"i100"
**@{--}; "Extra";"i001" **@{--}; "Extra";"100" **@{--};
"Extra";"001" **@{--};
"000";"i000" **@{--}; "001";"i001" **@{--}; "100";"i100" **@{--};
"110";"i110" **@{--}; "011";"i011" **@{--}; "111";"i111" **@{--};
\end{xy}
\]
\vskip -0.2cm
\caption{Schlegel diagram of a $4$-dimensional polytope with $10$
facets}\label{Oliver5}
\end{figure}

Each pair of adjacent permutations in a given class of $\Gamma$
specifies a CI statement. For instance, the four-element class
$\,45|1|23\,$ specifies the two CI statements
$ \ciaa{4}{5}{\emptyset}$ and $\ciaa{2}{3}{145}$,
while the $12$-element class $\,15|234\,$
specifies the seven CI statements
$$\,\ciaa{1}{5}{\emptyset}, \ciaa{2}{3}{15}, \ciaa{2}{3}{145}, \ciaa{2}{4}{15}, \ciaa{2}{4}{135}, \ciaa{3}{4}{15}, \ciaa{3}{4}{125}.$$
Altogether, we obtain $44$ CI statements $\,[\![ \cdot\, |\, \cdot
]\!] \,$ from the $14$ classes, and we identify the pre-convex rank
test $\Gamma$ with this set of $44$ CI statements. We now prove:

\begin{thm}
\label{thm:coarsest}
 $\Gamma$ is a coarsest convex rank test which is not submodular.
\end{thm}

\begin{proof}
To establish this theorem, we must prove the following three claims:
\begin{itemize}
\item $\Gamma$ is a convex rank test, i.e.~it satisfies the semigraphoid axioms ({\bf SG}).
\item There is no proper convex rank test which is coarser than $\Gamma$.
\item The convex rank test $\Gamma$ is not submodular.
\end{itemize}

\noindent We shall prove all three statements at once, by examining
the semigraphoid equations $({\bf SG}'')$. As in Section 2, the $44$
symbols in $\Gamma$ are denoted with double brackets $\,[\![ \,\cdot
\,|\,\cdot \,]\!]$, while the $36$ symbols not in $\Gamma$ are
denoted with brackets $\,[ \,\cdot \,|\,\cdot \,]$. With this distinction
between brackets, there are four symmetry types of semigraphoid equations that involve the 
$36$ positive unknowns $\,[ \,\cdot \,|\,\cdot \,]$. The full list is given in the Appendix:
$$ \begin{matrix}
{\rm Type \ I} &&& \cia{3}{5}{12} + \ciaa{3}{4}{125} & \,\, = \,\, & \cia{3}{4}{12} + \ciaa{3}{5}{124}  \\
{\rm Type \ II} &&& \cia{1}{5}{2} + \cia{1}{3}{25}        & \,\, = \,\, & \ciaa{1}{3}{2} + \cia{1}{5}{23}   \\
{\rm Type \ III} &&&   \cia{4}{5}{1} + \cia{2}{5}{14} & \,\,= \,\,& \cia{2}{5}{1} + \cia{4}{5}{12} \\
{\rm Type \ IV} &&&   \cia{1}{2}{5} + \ciaa{2}{3}{15} & \,\,= \,\,& \ciaa{2}{3}{5} + \cia{1}{2}{35}
   \end{matrix}  $$
  After setting the $44$ unknowns $\,[\![\,\cdot \,| \,\cdot \,]\!]\,$ to zero,
  we are left with $120$
equations in the $36$ strictly positive unknowns. For instance, the first 
three types give
$$ \begin{matrix}
{\rm Type \ I} &&& \cia{3}{5}{12}  & \,\, = \,\, & \cia{3}{4}{12}  \\
{\rm Type \ II} &&& \cia{1}{5}{2} + \cia{1}{3}{25} & \,\, = \,\, &  \cia{1}{5}{23}   \\
{\rm Type \ III} &&&   \cia{4}{5}{1} + \cia{2}{5}{14} & \,\,= \,\,& \cia{2}{5}{1} + \cia{4}{5}{12}
  \end{matrix}  $$

  The axiom $({\bf SG}'')$ merely requires that each of these equations is {\em individually} solvable.
  This is obviously the case. Hence $\Gamma$ is a semigraphoid.

The $78$ equations of Type I listed in the Appendix imply that all
$36$ positive unknowns must be equal. So, if another CI
statement is added to the semigraphoid $\Gamma$, then all others
must be added in order for ({\bf SG}) to remain valid.
  This proves our second claim that $\Gamma$ is a coarsest
convex rank test.

Given that the $36$ unknowns $\,[\,\cdot \, | \, \cdot \, ]$ must be
equal, the $12$ Type II equations imply that their common value is
zero, contradicting the requirement that they be positive. Hence the
$120$ orginal equations {\em altogether} have no non-negative real
solution that is consistent with $\Gamma$. This proves our third
claim that $\Gamma$ is not submodular.
\end{proof}

Every semigraphoid for $n=5$ corresponds to a $4$-dimensional fan.
Intersecting this fan with a sphere around the origin, we obtain a
polyhedral cell decomposition of the $3$-dimensional sphere. We do
not know whether each of these $3$-spheres can be realized as the
boundary of a $4$-dimensional polytope. However, using \cite[\S
5]{Zie}, every semigraphoid can be represented by a $3$-dimensional
diagram as in Figure 2.

For the specific semigraphoid $\Gamma$ of Theorem \ref{thm:coarsest},
the diagram in Figure 2 is indeed the boundary
of a $4$-polytope with f-vector $(14, 36, 32,
10)$.  The following coordinates for this polytope
 were found by a direct calculation, using
the techniques described in
\cite{Bokowski}. Each of the following ten row vectors
represents a facet of our polytope:
\begin{verbatim}
POINTS
 1      1/4         0       0       0        0  
 1       0          1       0       0        0  
 1       0          0       1       0        0  
 1       0          0       0       1        0  
 1       0          0       0       0        1  
 1     -1/4        1/4     1/4     5/4      1/4 
 1    280/893   -280/893  25/893    0      28/893
 1      1/57       1/57   -1/57   17/19     2/57
 1       1          1       0      -5        1  
 1      2/37      20/37    1/37   10/37    -2/37
 \end{verbatim}
 \vskip -0.4cm
For instance, the last row represents the facet-defining inequality
$$ \frac{2}{37} \cdot x_1 \,+\, \frac{20}{37} \cdot x_2 \,+\, \frac{1}{37} \cdot x_3
\,+\, \frac{10}{37} \cdot x_4 \,-\, \frac{2}{37} \cdot x_5 \,\,\, \leq \,\,\, 1. $$
Here, we are considering the vectors $(x_1,x_2,x_3,x_4,x_5)$ 
to be elements  in the quotient of $\mathbb{R}^5$ modulo the
one-dimensional linear subspace spanned by $(4,1,1,1,1)$.
Our format is that of the software {\tt Polymake} \cite{GJ}. 
If the above eleven lines are put in a file named
{\tt mypolytope} then the following command in {\tt Polymake}
will verify that this polytope does indeed have the 
combinatorial structure displayed in Figure 2:
\begin{verbatim}
 polymake mypolytope F_VECTOR VERTICES_IN_FACETS
 \end{verbatim}
 
 \vskip -0.4cm

The $10$ facets of our $4$-polytope correspond to the facets of the $5$-cube, and
they comprise all classes of permutations in $S_5$ in which the first or last
coordinate is fixed. The facets corresponding to permutations with
{\em $1$ or $5$} in the {\em first} coordinate have seven vertices,  twelve edges, and eight
2-faces.  The facets corresponding to permutations with {\em $2$,
$3$ or $4$ first} have seven vertices, $13$ edges, and eight
2-faces. The facets for {\em $1$ or $5$ last} are tetrahedra. The
 facets for {\em $2$, $3$ or $4$ last} are cubes in which one
edge has been contracted; they have seven vertices and $11$ edges.

\section{The semigraphoid semigroup is not normal} \label{sec:nonnormal}

Continuing to assume $n=5$, we now consider the linear map
$\mathcal{A}$ in the Introduction. It maps the free abelian group
$\mathbb{Z}^{80}$ spanned by the CI statements to the free abelian
group $\mathbb{Z}^{32}$ with basis $\,\{ e_K : K \subseteq [5] \}\,$
as specified in (\ref{Amap}). The matrix representing $\mathcal{A}$
has $32$ rows and $80$ columns; each column has four non-zero
entries: two $+1$'s and two $-1$'s. The rank of $\mathcal{A}$ is
$26$. The {\em semigraphoid semigroup}
 is $\, \mathcal{A}(\mathbb{N}^{80}) $, the non-negative integer span of
 the columns of this $32 \times 80$-matrix. This is
 a subsemigroup  of $\mathbb{Z}^{32}$. Equivalently,
the semigraphoid semigroup is the affine semigroup with $80$
generators and $120$ relations (given in the Appendix).  Note that
the polyhedral cone dual to the semigraphoid semigroup is the cone
of submodular functions.

In the language of \cite{Studeny2005Probabilistic}, the vectors in
$\mathbb{Z}^{32}$ are called {\em imsets}, the columns of
$\mathcal{A}$ are {\em elementary imsets}, and the elements of $
\mathcal{A}(\mathbb{N}^{80}) $ are {\em combinatorial imsets}. A
{\em structural imset} is a lattice point which lies in the
polyhedral cone spanned by the elementary imsets. Studen\' y's
question ({\bf Q2}) whether each structural imset is combinatorial
translates into the question whether the semigroup  $
\mathcal{A}(\mathbb{N}^{80}) $ is normal.

\begin{thm} \label{notnormal} The semigraphoid semigroup is not normal for $n=5$.
\end{thm}

\begin{proof}
Consider the following element in the free abelian group
$\mathbb{Z}^{80}$:
\begin{equation}
\label{exx1}
\begin{matrix}
& \cia{1}{5}{2}+\cia{1}{4}{3}+\cia{2}{3}{4}+\cia{2}{3}{5}+\cia{3}{4}{12} \\
& +\cia{2}{5}{13}
  +\cia{1}{2}{45}+\cia{1}{3}{45}+\cia{4}{5}{23} - \cia{2}{3}{45}.
 \end{matrix}
\end{equation}
The image of this element under the map $\mathcal{A} :
\mathbb{Z}^{80} \rightarrow \mathbb{Z}^{32}$ is the imset
\begin{equation}
\label{exx2}
\begin{matrix}
{\bf b} \,\, : = &
 -e_{2}-e_{3 }-e_{4 }-e_{5 }
-e_{23 }+e_{24 } +2e_{25 }+2e_{34 }+e_{35 }-e_{45 }
 +2e_{123 }
\\ & +e_{124 }-e_{125 }-e_{134 }+e_{ 135}+2e_{ 145}
-e_{1234 }-e_{1235 }-e_{1245 }-e_{1345 }.
\end{matrix}
\end{equation}
The imset ${\bf b}$ is structural because $2 \cdot {\bf b}$ is a
combinatorial imset. It is the image of
\begin{equation}
\label{exx3}
\begin{matrix}
 \cia{4}{5}{2}+\cia{4}{5}{3}+\cia{1}{3}{4}+\cia{1}{2}{5}+\cia{2}{5}{14}  +\cia{3}{4}{15}
&  \\
 +\cia{1}{4}{23}+\cia{1}{5}{23} +  \cia{1}{5}{2}+\cia{1}{4}{3} +\cia{2}{3}{4}
 & \\
 +\cia{2}{3}{5}+\cia{3}{4}{12}+\cia{2}{5}{13}+\cia{1}{2}{45}+\cia{1}{3}{45} & \in \,\, \mathbb{N}^{80}
 \end{matrix}
\end{equation}
under the linear map $\mathcal{A}$.

  Suppose that ${\bf b}$ were a combinatorial imset. Then there exists
   ${\bf x} \in \mathbb{N}^{80}$ such
that $\,\mathcal{A} \cdot {\bf x} = {\bf b}$. We write $\, {\bf x} \, =
\,\sum_i \cia{a_i}{b_i}{K_i} $, where we allow repetition in the
sum. In any elementary imset, the basis vector $e_\emptyset$ occurs
with coefficient $-1$ or $0$,
 and the basis vector $e_{12345}$ occurs with coefficient $-1$ or $0$.
However, neither $\,e_\emptyset\,$ nor $\,e_{12345}\,$ appears in the
imset ${\bf b}$, so we conclude that $| K_i | = 1$ or $|K_i| = 2$
for all terms $\,\cia{a_i}{b_i}{K_i}\,$ in the representation of
${\bf x}$. The first four terms $\, -e_{2}-e_{3 }-e_{4 }-e_{5 }\,$ in
$\,{\bf b}\,$
 imply that ${\bf x}$ has precisely four terms $\,\cia{a_i}{b_i}{K_i}\,$
with $|K_i|= 1$, and the terms $\,-e_{1234 }-e_{1235 }-e_{1245
}-e_{1345 }\,$ imply that $\,{\bf x}\,$ has precisely four terms
with $|K_i|= 2$.

Each of the eight terms in $\,{\bf x}\,$ evaluates to an alternating
sum of $4$ terms under the map $\mathcal{A}$. Some
cancellation occurs among the resulting $32$ terms. Prior to that
cancellation, our imset  had been written as the sum of two subsums,
$\, {\bf b} \, = \, \mathcal{A} \cdot {\bf x}\,= $
$$
 -e_{ 2}-e_{3}-e_{4}-e_{5}
+ e_{24} +2e_{25}+2e_{34}+e_{35}+ e_{A_1} +e_{A_2} -e_{125}-e_{134}
- e_{B_1} - e_{B_2}  $$
$$
 -e_{23} - e_{45} - e_{A_1}  - e_{A_2}
 + 2e_{123} + e_{124} +e_{ 135} + 2e_{ 145} + e_{B_1} + e_{B_2}
-e_{1234}-e_{1235}-e_{1245}-e_{1345},
$$
\vskip 0.2cm \noindent where $\,|A_1| = |A_2| = 2 \,$ and
 $\,|B_1| = |B_2| = 3 $. The first line
is the sum of the four elementary imsets $\,
\mathcal{A}(\cia{a_i}{b_i}{K_i} )\,$ with $|K_i|=1$, and the second
line is the sum of the four elementary imsets with $|K_i|=2$. A
contradiction will arise when we try to determine the unknown pairs
$A_1$ and $A_2$. The term $-e_{125}$ in the first line must come from
$K_i = \{2\}$ or $K_i = \{5\}$. This implies that either $\{1,2\}$
or $\{1,5\}$ is in $\,A_* = \{A_1,A_2\}$. Similarly, the term
$-e_{134}$ shows that  either $\{1,3\}$ or $\{1,4\}$ is in $A_*$. Now
consider the second line. The presence of the term $2 e_{123}$
implies that $\{1,2 \}$ or $\{1,3\}$ is in $A_*$, and the term  $2
e_{145}$  implies that $\{1,4 \}$ or $\{1,5\}$ is in $A_*$. The term
$e_{124}$ shows that $\{1,2 \}$, $\{ 1,4 \}$, or $\{2,4\}$ is in
$A_*$, and, finally, the term $e_{135}$ shows that $\{1,3 \}$, $\{
1,5\}$, or $\{3,5\}$ is in $A_*$. However, no such pair of pairs
$A_*$ satisfies these six restrictions. This proves that ${\bf b}$
is not a combinatorial imset.
\end{proof}

The main point of the above proof was to show that the linear system
$\, \mathcal{A} \cdot {\bf x} \, = \,{\bf b}\,$ has no solution
with non-negative {\em integer} coordinates. This can also be
verified automatically  using integer programming software. In fact,
using such software we found that $\, \mathcal{A} \cdot {\bf x} \, =
\,{\bf b}\,$ has only one solution with non-negative {\em real}
coordinates, namely, that unique  solution $\,{\bf x} \in
(\mathbb{R}_{\geq 0})^{80}\,$ is the expression in (\ref{exx3})
scaled by $1/2$.

The reader might now inquire how the imset ${\bf
b}$ was found. There are several algorithms that test whether a
given affine semigroup is normal, including one recently proposed by
Takemura, Yoshida and the first author \cite{HTY}, and the method of
Bruns and Koch \cite{BK} which is implemented in their software {\tt
normaliz}.

Our original attempts to apply these methods directly to the $32 \times
80$-matrix $\mathcal{A}$ were unsuccessful.
Instead we succeeded by partially computing a
   {\em Markov basis} for the matrix
$\mathcal{A}$ using the software {\tt 4ti2} \cite{HHM}. The imset
{\bf b} was found by inspecting the partial results produced by {\tt
4ti2}. We explain the details in the next section.

\section{Computations in toric algebra}

Let $\mathbb{Q}[{\rm CI}_n]$ denote the polynomial ring over the
field of rational numbers $\mathbb{Q}$ generated by the symbols $\,
[ i \indep j \, |\, K]$. Thus $\mathbb{Q}[{\rm CI}_n]$ is a
polynomial ring in $\gamma_n$ unknowns, one for each $2$-face of the
$n$-cube $C_n$. We write $\,\prod {\rm CI}_n\,$ for the product of
all the unknowns. We define the {\em semigraphoid ideal} to be the
ideal $I_{{\bf SG}}$ generated by the binomials in (${\bf SG}'''$).
Thus the generators of $I_{{\bf SG}}$ represent the semigraphoid
axioms.  Following \cite[\S 7]{MS}, we introduce the {\em toric
ideal} $I_{\mathcal{A}}$ which is obtained from $I_{{\bf SG}}$ by
saturation:
\begin{equation}
\label{saturate}
 I_{\mathcal{A}} \quad := \quad \bigl( I_{{\bf SG}} \,: \, (\prod {\rm CI}_n)^\infty \bigr).
 \end{equation}
The binomials in  $I_{\mathcal{A}}$ represent the vectors in the
kernel of the linear map $\,\mathcal{A} : \mathbb{Z}^{\gamma_n}
\rightarrow \mathbb{Z}^{2^n}$. A minimal set of binomials which
generates $I_{\mathcal{A}}$ is said to be a {\em Markov basis} for
the matrix $\mathcal{A}$. See \cite{DS} for a discussion of Markov
bases in the context of statistics.

Let us illustrate these concepts for $n = 3$. The polynomial ring
$\mathbb{Q}[{\rm CI}_3]$ has six unknowns, one for each facet of the
$3$-cube. They are the entries of the $2 \times 3$-matrix
\begin{equation}
\label{twobythree}
\begin{pmatrix}
\cia{1}{2}{\emptyset} & \cia{1}{3}{\emptyset} &
\cia{2}{3}{\emptyset} \\
\cia{1}{2}{3} & \cia{1}{3}{2} & \cia{2}{3}{1}
\end{pmatrix}.
\end{equation}
The semigraphoid ideal $I_{{\bf SG}}$ is generated by the three $2
\times 2$-minors of the matrix (\ref{twobythree}). This is a prime
ideal of codimension $2$ and degree $3$, and hence we have $\,
I_{{\bf SG}} \, = \, I_{\mathcal{A}}$. Here the Markov basis for
$\mathcal{A}$ consists precisely of the three semigraphoid axioms.

We next consider the case $n=4$. The polynomial ring
$\mathbb{Q}[{\rm CI}_4]$ has $24$ unknowns, one for each $2$-face of
the $4$-cube. They are the entries of eight $2 \times 3$-matrices as
in (\ref{twobythree}), one for each of the eight facets of the
$4$-cube. Thus the semigraphoid ideal $I_{{\bf SG}}$ is generated by
$24$ quadrics, one for each of the $24$ axioms (${\bf SG}''$) in the
list given in Section 2. For instance, the last axiom in that list
translates into the quadratic binomial $\,[1
\!\!\perp\!\!\!\perp\!\!2|4] \cdot [2
\!\!\perp\!\!\!\perp\!\!3|14]\, - \,[2
\!\!\perp\!\!\!\perp\!\!3|4] \cdot [1 \!\!\perp\!\!\!\perp\!\!2|34]
$, which is one of the $24$ generators of $I_{{\bf SG}}$. Using the
software {\tt Macaulay2} \cite{M2} we derived the following result:

\begin{prop}
\label{primedec} The semigraphoid ideal $I_{{\bf SG}}$ is a radical
ideal which is the intersection of the toric ideal $I_{\mathcal{A}}$
and $17$ additional associated monomial prime ideals.
\end{prop}

Before discussing this prime decomposition in detail, let us
make a few general remarks. We wish to argue that toric algebra and
algebraic geometry provide useful algorithmic tools for the research
directions presented in \cite{Studeny2005Probabilistic}. For any
ideal $I $ of $\mathbb{Q}[{\rm CI}_n]$ and any subset
 $\Omega$ of the complex affine space $\mathbb{C}^{\gamma_n}$,
 the {\em variety} $\,V_\Omega(I)\,$ is defined as
 the set of all vectors in $\Omega$
 which are common zeros of all the polynomials in $I$.
 Then $V_{\C}(I_{{\bf SG}})$ is a complex variety, reducible for $n \geq 4$, one
 of whose irreducible components is the complex toric variety
  $\,V_{\C} (I_{\mathcal{A}})$. Inside this toric variety are
  the real toric variety $\,V_{\R} (I_{\mathcal{A}})$. Its
  non-negative part $V_{{\R}_{\geq 0}}(I_\mathcal{A})$
  is homeomorphic to the cone spanned by the elementary imsets.
  Our next result shows that
  the semigraphoids are precisely the points on these varieties
whose coordinates are $0$ or~$1$.

\begin{thm}
\label{rmkkk} The semigraphoids on $[n]$ are in bijection with the
points in $V_{\{0,1\}}(I_{{\bf SG}})$. The submodular semigraphoids
on $[n]$ are in bijection with the points in
$\,V_{\{0,1\}}(I_{\mathcal{A}})$.
\end{thm}

\begin{proof}
We replace the additive semiring $\{\szero,\splus\}$ with the
multiplicative semiring $\{1,0\}$. This translates from the
 additive notation (${\bf SG}''$) to the
multiplicative notation (${\bf SG}'''$). With this translation,
the first statement in Theorem \ref{rmkkk} is obvious.

The second statement is less obvious and is based on the geometry
of toric varieties. Specifically, we shall use the characterization of
{\em facial} index sets which is developed in \cite{GMS}. If we 
consider our specific $2^n \times \gamma_n$-matrix $\mathcal{A}$ then the role of the
set $\{1,\ldots,m\}$ in \cite{GMS} is played by the
set of CI statements, and a subset of CI statements is facial for $\mathcal{A}$ if and only if
it is submodular semigraphoid. With this observation,
our second assertion follows from 
Lemma A.2 in the Appendix of \cite{GMS}.
\end{proof}

 Using Theorem \ref{rmkkk}, we can study semigraphoids by studying
 the zero-dimensional ideals obtained by adding $\, \langle\, x^2-x \,: \, x \in {\rm CI}_n \,\rangle \,$
 to the ideal $I_{{\bf SG}}$ or $I_{\mathcal{A}}$. For instance, with
  the command {\tt degree} in {\tt Macaulay2} \cite{M2}, it takes only a few seconds to compute
\begin{equation}
\label{QuickCount} \# V_{\{0,1\}}(I_{{\bf SG}}) \, = \, 26424 \qquad
\hbox{and} \qquad \# V_{\{0,1\}}(I_{\mathcal{A}}) \, = \,  22108.
 \end{equation}
The difference between these numbers is explained geometrically
by the prime decomposition in Proposition \ref{primedec}, which we
shall now describe in explicit terms.

The $17$ associated monomial primes of $I_{{\bf SG}}$ come in three
symmetry classes. First there are two primes of codimension $12$. A
representative is the ideal
$$ \left\langle
\begin{matrix}
\cia{1}{2}{\emptyset}, \, \cia{1}{3}{\emptyset},\,
\cia{1}{4}{\emptyset}, \, \cia{2}{3}{\emptyset},\,
\cia{2}{4}{\emptyset}, \,
 \cia{3}{4}{\emptyset},
\\ \,\cia{3}{4}{12}, \cia{2}{4}{13},
\cia{2}{3}{14}, \cia{1}{4}{23},
    \cia{1}{3}{24}, \cia{1}{2}{34} \,
\end{matrix} \right\rangle . $$
The semigraphoid ideal $I_{{\bf SG}}$  has $12$ associated primes of
codimension $15$, such as
$$  \left\langle
\begin{matrix} \,
\cia{1}{2}{\emptyset}, \cia{1}{3}{\emptyset}, \cia{1}{4}{\emptyset},
\cia{3}{4}{\emptyset}, \cia{1}{3}{2},
 \cia{1}{4}{2}, \cia{3}{4}{2}, \cia{1}{2}{3}, \\
  \cia{2}{4}{3}, \cia{1}{2}{4},
    \cia{2}{3}{4}, \cia{3}{4}{12}, \cia{2}{4}{13}, \cia{2}{3}{14}, \cia{1}{2}{34}
\end{matrix} \right\rangle . $$
Next, $I_{{\bf SG}}$ has three associated primes of codimension
$16$. A representative is
$$ \left\langle
\begin{matrix} \,
\cia{1}{2}{\emptyset}, \cia{1}{3}{\emptyset}, \cia{2}{4}{\emptyset},
\cia{3}{4}{\emptyset}, \cia{2}{4}{1}, \cia{3}{4}{1},
 \cia{1}{3}{2}, \cia{3}{4}{2}, \\
 \cia{1}{2}{3}, \cia{2}{4}{3},
    \cia{1}{2}{4}, \cia{1}{3}{4}, \cia{3}{4}{12}, \cia{2}{4}{13}, \cia{1}{3}{24}, \cia{1}{2}{34}
\end{matrix} \right\rangle. $$
Each of the $ 4316$ non-submodular semigraphoids is a
$\{0,1\}$-valued point not in $V(I_{\mathcal{A}})$ but in one of the
$17$ coordinate subspaces corresponding to these primes.

Finally, the last associated prime of $I_{{\bf SG}}$ is the toric
ideal $I_{\mathcal{A}}$. This ideal has codimension $13$ and degree
$396$. Its minimal generating set consists of $52$ binomials.
Besides the
 $24$ quadrics (the axioms), the Markov basis of $\mathcal{A}$ contains
four cubics
$$ \begin{matrix}
   \cia{2}{3}{1} \cdot \cia{3}{4}{2} \cdot \cia{1}{3}{4} - \cia{3}{4}{1} \cdot \cia{1}{3}{2} \cdot \cia{2}{3}{4}, \\
  \cia{2}{3}{1} \cdot \cia{2}{4}{3} \cdot \cia{1}{2}{4} - \cia{2}{4}{1} \cdot \cia{1}{2}{3} \cdot \cia{2}{3}{4}, \\
  \cia{1}{3}{2} \cdot \cia{1}{4}{3} \cdot \cia{1}{2}{4} - \cia{1}{4}{2} \cdot \cia{1}{2}{3} \cdot \cia{1}{3}{4}, \\
  \cia{2}{4}{1} \cdot \cia{3}{4}{2} \cdot \cia{1}{4}{3} - \cia{3}{4}{1} \cdot \cia{1}{4}{2} \cdot \cia{2}{4}{3},
    \end{matrix}
 $$
and $24$ quartics such as
$$
  \cia{1}{2}{\emptyset} \cdot \cia{3}{4}{\emptyset} \cdot \cia{2}{4}{13} \cdot \cia{1}{3}{24}
 \,  - \, \cia{1}{3}{\emptyset} \cdot \cia{2}{4}{\emptyset} \cdot \cia{3}{4}{12} \cdot \cia{1}{2}{34}.
 $$

We now come to case $n=5$. It will be a challenge for future commutative
algebra software to compute a primary decomposition of the semigraphoid ideal $\,I_{\bf SG}\,$  
for $n=5$. At present we do not know even whether $\,I_{\bf SG}\,$
is radical. Let us therefore focus on the main component of this
ideal, namely, the toric ideal $I_{\mathcal{A}}$. Here our main goal
is to compute its minimal generators, that is, the Markov basis of
$\mathcal{A}$. We attacked this problem using the software {\tt
4ti2} \cite{HHM}, and we now discuss the results.

First, we started a Markov basis computation for the toric ideal
$I_{\mathcal{A}}$ using the function {\tt markov} of {\tt 4ti2}, but this computation turned out to be non-trivial. 
In the hope that a counterexample would not involve all $80$
variables, we set several variables to $0$ and tried to compute the
Markov basis of smaller ideals that are contained in
$I_{\mathcal{A}}$. For the one-day computation that finally produced
a counterexample, we set the first $18$ formal symbols to zero and found
the Markov basis move
\[
{\bf g} \,\,\, := \,\,\, \bigl(\,\mathbf{\alpha}+2 \cdot \cia{2}{3}{45}\,\bigr)\,-\, \bigl(\,\mathbf{\beta}+2 \cdot \cia{4}{5}{23} \,\bigr)\,\,
\in \,\, \mathbb{N}^{80}, \qquad \qquad \hbox{where} \qquad
\]
{\footnotesize
\[
\mathbf{\alpha} =
\cia{4}{5}{2}+\cia{4}{5}{3}+\cia{1}{3}{4}+\cia{1}{2}{5}+\cia{2}{5}{14}
+\cia{3}{4}{15} +\cia{1}{4}{23}+\cia{1}{5}{23},
\]
\[
\mathbf{\beta} = \cia{1}{5}{2}+\cia{1}{4}{3} +\cia{2}{3}{4}
+\cia{2}{3}{5}+\cia{3}{4}{12}+\cia{2}{5}{13}+\cia{1}{2}{45}+\cia{1}{3}{45}.
\]
}
This lattice vector corresponds to a binomial $\,{\bf x}^{{\bf
g}^+}-\,{\bf x}^{{\bf g}^-}\,$ which is in the toric ideal $\, I_{\mathcal{A}}\,$ 
and has the property that both of its monomials are not square-free.
 We then verified that ${\bf x}^{{\bf
g}^+}-\,{\bf x}^{{\bf g}^-}$ is not only indispensable for the
smaller ideal (with $18$ variables set to zero) but also
indispensable for $I_{\mathcal{A}}$. 
Recall (e.g.~from \cite{ATY}) that a binomial ${\bf
x}^{{\bf g}^+}-\,{\bf x}^{{\bf g}^-}\,$ in the toric ideal $\, I_{\mathcal{A}}\,$ is called
{\em indispensable} if
\[
\{{\bf z}\in\,\,\mathbb{N}^{80}:\mathcal{A}\cdot{\bf
z}=\mathcal{A}\cdot{\bf g}^+\} \quad = \quad \{{\bf g}^+,{\bf g}^-\}.
\]
This means that the Markov move ${\bf g}$ corresponds to a $2$-element
fiber given by the right-hand side and consequently, ${\bf g}$ must
belong to {\em every} Markov basis of $I_{\mathcal{A}}$. In
order to check this condition for our given move ${\bf g}$, we
computed the minimal Hilbert basis (that is, the $\leq$-minimal
integer solutions) of the cone
\[
\{({\bf z}, u)\in\,\,\mathbb{R}^{81}:\mathcal{A}\cdot{\bf
z}-(\mathcal{A}\cdot{\bf g}^+)\cdot u =0, ({\bf z},u )\geq
0\}.
\]
This was done using the function {\tt hilbert} of {\tt 4ti2} which
produced precisely the two expected elements $(\mathbf{g}^+,1)$ and $(\mathbf{g}^-,1)$
within a few seconds. 

From our special Markov move
$\, {\bf g} =(\mathbf{\alpha}+2 \cdot \cia{2}{3}{45})-(\mathbf{\beta}+2  \cdot \cia{4}{5}{23})$, we then
constructed the imset ${\bf b}$ presented in Section 4.
We first checked that $\mathbf{b}$ was not a combinatorial imset by
showing that $\mathcal{A}{\bf x} = {\bf b}$ has no solutions with 
non-negative integer coordinates.
Using the functions \ {\tt hilbert} \ and 
\ {\tt rays} \ of the program
\ {\tt 4ti2}, we computed the Hilbert basis and the extreme rays of
the cone
\[
\bigl\{({\bf z}, u)\in\,\,\mathbb{R}^{81}:\mathcal{A}\cdot{\bf z} =
{\bf b}\cdot u \,\,\hbox{and} \,\, ({\bf z}, u)\geq 0 \bigr\}.
\]
Both computations quickly finished. They showed that this cone 
has dimension one and is generated by the single vector $({\bf \alpha}+{\bf \beta},2)$. Consequently,
the only non-negative real solution to $\, \mathcal{A} \cdot {\bf x}
\, = \,{\bf b}\,$ is $({\bf \alpha}+{\bf \beta})/2$, which is not an integer solution.

\vskip .1cm

We are currently in the process of computing the complete  minimal Markov basis
of the toric ideal for semigraphoids with $n=5$.
That Markov basis has well over a million elements. Yet,
we are convinced that {\tt 4ti2} will succeed.
The completion of that Markov basis
will represent a computational milestone in toric algebra.

\medskip

\section*{Acknowledgments}

Jason Morton and Bernd Sturmfels were supported by the
DARPA {\em Fundamental Laws of Biology} program,
and Bernd Sturmfels was also supported by the NSF.
Anne Shiu was supported by a Lucent Technologies Bell Labs Graduate Research Fellowship.  Oliver Wienand was supported by the Wipprecht foundation.

\bigskip

\noindent {\bf Authors' addresses:}

\medskip

\noindent

\noindent Raymond Hemmecke,
Fakult\"at f\"ur Mathematik,
Otto-von-Guericke-Universit\"at \break Magdeburg,
39106 Magdeburg, Germany,
{\tt raymond@hemmecke.de}

\medskip

\noindent Jason Morton, Anne Shiu and Bernd Sturmfels,
 Department of Mathematics, \break
University of California,  Berkeley, CA 94720, USA, \hfill \break
{\tt [mortonj, annejls, bernd]@math.berkeley.edu}

\medskip

\noindent Oliver Wienand, 
Fachbereich Mathematik,
Universit\"at Kaiserslautern,
67653 Kaiserslautern, Germany,
{\tt wienand@rhrk.uni-kl.de}

\vfill \eject

\section{Appendix: The 120 semigraphoid axioms}

Here is the list of all $120 $ semigraphoid axiom for
$n=5$, grouped into triples according to which $3$-face of the
$5$-cube they come from.  The two types of brackets specify the
non-submodular coarsest semigraphoid $\Gamma$ 
which was discussed in Section
\ref{sec:coarsest}.
{\tiny
$$
\begin{array}{rclcrcl}
  \cia{3}{5}{12} + \ciaa{3}{4}{125} & = & \cia{3}{4}{12} + \ciaa{3}{5}{124} &
& \cia{2}{5}{13} + \ciaa{2}{4}{135} & = & \cia{2}{4}{13} + \ciaa{2}{5}{134} \\
  \cia{4}{5}{12} + \ciaa{3}{4}{125} & = & \cia{3}{4}{12} + \ciaa{4}{5}{123} &
& \cia{4}{5}{13} + \ciaa{2}{4}{135} & = & \cia{2}{4}{13} + \ciaa{4}{5}{123} \\
  \cia{4}{5}{12} + \ciaa{3}{5}{124} & = & \cia{3}{5}{12} + \ciaa{4}{5}{123} &
& \cia{4}{5}{13} + \ciaa{2}{5}{134} & = & \cia{2}{5}{13} + \ciaa{4}{5}{123} \\
\\
  \cia{2}{3}{14} + \ciaa{2}{5}{134} & = & \cia{2}{5}{14} + \ciaa{2}{3}{145} &
& \ciaa{2}{4}{15} + \ciaa{2}{3}{145} & = & \ciaa{2}{3}{15} + \ciaa{2}{4}{135}\\
  \cia{2}{5}{14} + \ciaa{3}{5}{124} & = & \cia{3}{5}{14} + \ciaa{2}{5}{134} &
& \ciaa{3}{4}{15} + \ciaa{2}{3}{145} & = & \ciaa{2}{3}{15} + \ciaa{3}{4}{125}\\
  \cia{3}{5}{14} + \ciaa{2}{3}{145} & = & \cia{2}{3}{14} + \ciaa{3}{5}{124} &
& \ciaa{3}{4}{15} + \ciaa{2}{4}{135} & = & \ciaa{2}{4}{15} + \ciaa{3}{4}{125}\\
\\
  \cia{1}{5}{23} + \ciaa{1}{4}{235} & = & \cia{1}{4}{23} + \ciaa{1}{5}{234}   &
& \cia{1}{3}{24} + \ciaa{3}{5}{124} & = & \cia{3}{5}{24} + \ciaa{1}{3}{245}\\
   \cia{4}{5}{23} + \ciaa{1}{4}{235} & = & \cia{1}{4}{23} + \ciaa{4}{5}{123}  &
& \cia{1}{5}{24} + \ciaa{1}{3}{245} & = & \cia{1}{3}{24} + \ciaa{1}{5}{234}\\
   \cia{4}{5}{23} + \ciaa{1}{5}{234} & = & \cia{1}{5}{23} + \ciaa{4}{5}{123}&
& \cia{1}{5}{24} + \ciaa{3}{5}{124} & = & \cia{3}{5}{24} + \ciaa{1}{5}{234}\\
\\
 \cia{1}{3}{25} + \ciaa{1}{4}{235} & = & \cia{1}{4}{25} + \ciaa{1}{3}{245}     &
& \cia{1}{2}{34} + \ciaa{2}{5}{134} & = & \cia{2}{5}{34} + \ciaa{1}{2}{345}\\
 \cia{1}{3}{25} + \ciaa{3}{4}{125} & = & \cia{3}{4}{25} + \ciaa{1}{3}{245}    &
& \cia{1}{5}{34} + \ciaa{1}{2}{345} & = & \cia{1}{2}{34} + \ciaa{1}{5}{234}\\
 \cia{3}{4}{25} + \ciaa{1}{4}{235} & = & \cia{1}{4}{25} + \ciaa{3}{4}{125}&
& \cia{1}{5}{34} + \ciaa{2}{5}{134} & = & \cia{2}{5}{34} + \ciaa{1}{5}{234}\\
\\
  \cia{1}{2}{35} + \ciaa{1}{4}{235} & = & \cia{1}{4}{35} + \ciaa{1}{2}{345}   &
& \cia{1}{2}{45} + \ciaa{1}{3}{245} & = & \cia{1}{3}{45} + \ciaa{1}{2}{345}\\
  \cia{1}{2}{35} + \ciaa{2}{4}{135} & = & \cia{2}{4}{35} + \ciaa{1}{2}{345}   &
& \cia{1}{3}{45} + \ciaa{2}{3}{145} & = & \cia{2}{3}{45} + \ciaa{1}{3}{245}\\
  \cia{1}{4}{35} + \ciaa{2}{4}{135} & = & \cia{2}{4}{35} + \ciaa{1}{4}{235}&
& \cia{2}{3}{45} + \ciaa{1}{2}{345} & = & \cia{1}{2}{45} + \ciaa{2}{3}{145}\\
\\
  \ciaa{2}{4}{1} + \cia{3}{4}{12} & = & \ciaa{3}{4}{1} + \cia{2}{4}{13}    &
& \ciaa{2}{3}{1} + \cia{3}{5}{12} & = & \cia{3}{5}{1} + \ciaa{2}{3}{15}\\
  \ciaa{2}{4}{1} + \cia{2}{3}{14} & = & \ciaa{2}{3}{1} + \cia{2}{4}{13}    &
& \ciaa{2}{3}{1} + \cia{2}{5}{13} & = & \cia{2}{5}{1} + \ciaa{2}{3}{15}\\
  \ciaa{3}{4}{1} + \cia{2}{3}{14} & = & \ciaa{2}{3}{1} + \cia{3}{4}{12}&
& \cia{3}{5}{1} + \cia{2}{5}{13} & = & \cia{2}{5}{1} + \cia{3}{5}{12}\\
\\
  \ciaa{2}{4}{1} + \cia{4}{5}{12} & = & \cia{4}{5}{1} + \ciaa{2}{4}{15}    &
& \cia{3}{5}{1} + \cia{4}{5}{13} & = & \cia{4}{5}{1} + \cia{3}{5}{14}\\
  \cia{2}{5}{1} + \ciaa{2}{4}{15} & = & \ciaa{2}{4}{1} + \cia{2}{5}{14}    &
& \cia{3}{5}{1} + \ciaa{3}{4}{15} & = & \ciaa{3}{4}{1} + \cia{3}{5}{14}\\
  \cia{4}{5}{1} + \cia{2}{5}{14} & = & \cia{2}{5}{1} + \cia{4}{5}{12}&
& \cia{4}{5}{1} + \ciaa{3}{4}{15} & = & \ciaa{3}{4}{1} + \cia{4}{5}{13}\\
\\
   \ciaa{1}{3}{2} + \cia{1}{4}{23} & = & \ciaa{1}{4}{2} + \cia{1}{3}{24}    &
& \cia{1}{5}{2} + \cia{1}{3}{25} & = & \ciaa{1}{3}{2} + \cia{1}{5}{23}\\
   \ciaa{1}{4}{2} + \cia{3}{4}{12} & = & \ciaa{3}{4}{2} + \cia{1}{4}{23}    &
& \ciaa{3}{5}{2} + \cia{1}{5}{23} & = & \cia{1}{5}{2} + \cia{3}{5}{12}\\
   \ciaa{3}{4}{2} + \cia{1}{3}{24} & = & \ciaa{1}{3}{2} + \cia{3}{4}{12}&
& \ciaa{3}{5}{2} + \cia{1}{3}{25} & = & \ciaa{1}{3}{2} + \cia{3}{5}{12}\\
\\
  \ciaa{1}{4}{2} + \cia{4}{5}{12} & = & \ciaa{4}{5}{2} + \cia{1}{4}{25}    &
& \ciaa{3}{5}{2} + \cia{4}{5}{23} & = & \ciaa{4}{5}{2} + \cia{3}{5}{24}\\
  \ciaa{1}{4}{2} + \cia{1}{5}{24} & = & \cia{1}{5}{2} + \cia{1}{4}{25}    &
& \ciaa{3}{5}{2} + \cia{3}{4}{25} & = & \ciaa{3}{4}{2} + \cia{3}{5}{24}\\
  \cia{1}{5}{2} + \cia{4}{5}{12} & = & \ciaa{4}{5}{2} + \cia{1}{5}{24}&
& \ciaa{4}{5}{2} + \cia{3}{4}{25} & = & \ciaa{3}{4}{2} + \cia{4}{5}{23}\\
\\
  \ciaa{1}{4}{3} + \cia{2}{4}{13} & = & \ciaa{2}{4}{3} + \cia{1}{4}{23}    &
& \ciaa{1}{2}{3} + \cia{2}{5}{13} & = & \ciaa{2}{5}{3} + \cia{1}{2}{35}\\
  \ciaa{1}{4}{3} + \cia{1}{2}{34} & = & \ciaa{1}{2}{3} + \cia{1}{4}{23}    &
& \cia{1}{5}{3} + \cia{1}{2}{35} & = & \ciaa{1}{2}{3} + \cia{1}{5}{23}\\
  \ciaa{2}{4}{3} + \cia{1}{2}{34} & = & \ciaa{1}{2}{3} + \cia{2}{4}{13}&
& \ciaa{2}{5}{3} + \cia{1}{5}{23} & = & \cia{1}{5}{3} + \cia{2}{5}{13}\\
\\
  \ciaa{1}{4}{3} + \cia{4}{5}{13} & = & \ciaa{4}{5}{3} + \cia{1}{4}{35}    &
& \ciaa{2}{4}{3} + \cia{4}{5}{23} & = & \ciaa{4}{5}{3} + \cia{2}{4}{35}\\
  \ciaa{1}{4}{3} + \cia{1}{5}{34} & = & \cia{1}{5}{3} + \cia{1}{4}{35}    &
& \ciaa{2}{5}{3} + \cia{2}{4}{35} & = & \ciaa{2}{4}{3} + \cia{2}{5}{34}\\
  \ciaa{4}{5}{3} + \cia{1}{5}{34} & = & \cia{1}{5}{3} + \cia{4}{5}{13}&
& \ciaa{4}{5}{3} + \cia{2}{5}{34} & = & \ciaa{2}{5}{3} + \cia{4}{5}{23}\\
\\
  \ciaa{1}{2}{4} + \cia{2}{3}{14} & = & \ciaa{2}{3}{4} + \cia{1}{2}{34}    &
& \ciaa{1}{2}{4} + \cia{2}{5}{14} & = & \ciaa{2}{5}{4} + \cia{1}{2}{45}\\
  \ciaa{1}{2}{4} + \cia{1}{3}{24} & = & \ciaa{1}{3}{4} + \cia{1}{2}{34}    &
& \ciaa{1}{2}{4} + \cia{1}{5}{24} & = & \cia{1}{5}{4} + \cia{1}{2}{45}\\
  \ciaa{1}{3}{4} + \cia{2}{3}{14} & = & \ciaa{2}{3}{4} + \cia{1}{3}{24}&
& \cia{1}{5}{4} + \cia{2}{5}{14} & = & \ciaa{2}{5}{4} + \cia{1}{5}{24}\\
\\
  \ciaa{1}{3}{4} + \cia{1}{5}{34} & = & \cia{1}{5}{4} + \cia{1}{3}{45}    &
& \ciaa{2}{3}{4} + \cia{3}{5}{24} & = & \ciaa{3}{5}{4} + \cia{2}{3}{45}\\
  \ciaa{3}{5}{4} + \cia{1}{5}{34} & = & \cia{1}{5}{4} + \cia{3}{5}{14}    &
& \ciaa{2}{3}{4} + \cia{2}{5}{34} & = & \ciaa{2}{5}{4} + \cia{2}{3}{45}\\
  \ciaa{3}{5}{4} + \cia{1}{3}{45} & = & \ciaa{1}{3}{4} + \cia{3}{5}{14}&
& \ciaa{2}{5}{4} + \cia{3}{5}{24} & = & \ciaa{3}{5}{4} + \cia{2}{5}{34}\\
\\
  \cia{1}{2}{5} + \ciaa{2}{3}{15} & = & \ciaa{2}{3}{5} + \cia{1}{2}{35}    &
& \cia{1}{2}{5} + \ciaa{2}{4}{15} & = & \ciaa{2}{4}{5} + \cia{1}{2}{45}\\
  \cia{1}{2}{5} + \cia{1}{3}{25} & = & \cia{1}{3}{5} + \cia{1}{2}{35}    &
& \cia{1}{2}{5} + \cia{1}{4}{25} & = & \cia{1}{4}{5} + \cia{1}{2}{45}\\
  \cia{1}{3}{5} + \ciaa{2}{3}{15} & = & \ciaa{2}{3}{5} + \cia{1}{3}{25}&
& \cia{1}{4}{5} + \ciaa{2}{4}{15} & = & \ciaa{2}{4}{5} + \cia{1}{4}{25}\\
\\
  \cia{1}{3}{5} + \ciaa{3}{4}{15} & = & \ciaa{3}{4}{5} + \cia{1}{3}{45}    &
& \ciaa{2}{3}{5} + \cia{2}{4}{35} & = & \ciaa{2}{4}{5} + \cia{2}{3}{45}\\
  \cia{1}{3}{5} + \cia{1}{4}{35} & = & \cia{1}{4}{5} + \cia{1}{3}{45}    &
& \ciaa{2}{4}{5} + \cia{3}{4}{25} & = & \ciaa{3}{4}{5} + \cia{2}{4}{35}\\
  \ciaa{3}{4}{5} + \cia{1}{4}{35} & = & \cia{1}{4}{5} + \ciaa{3}{4}{15}&
& \ciaa{3}{4}{5} + \cia{2}{3}{45} & = & \ciaa{2}{3}{5} + \cia{3}{4}{25}\\
\\
  \ciaa{1}{2}{} + \ciaa{2}{3}{1} & = & \ciaa{2}{3}{} + \ciaa{1}{2}{3}    &
& \ciaa{1}{2}{} + \ciaa{2}{4}{1} & = & \ciaa{2}{4}{} + \ciaa{1}{2}{4}\\
  \ciaa{1}{3}{} + \ciaa{1}{2}{3} & = & \ciaa{1}{2}{} + \ciaa{1}{3}{2}    &
& \ciaa{1}{2}{} + \ciaa{1}{4}{2} & = & \ciaa{1}{4}{} + \ciaa{1}{2}{4}\\
  \ciaa{2}{3}{} + \ciaa{1}{3}{2} & = & \ciaa{1}{3}{} + \ciaa{2}{3}{1}&
& \ciaa{1}{4}{} + \ciaa{2}{4}{1} & = & \ciaa{2}{4}{} + \ciaa{1}{4}{2}\\
\\
  \ciaa{1}{2}{} + \cia{2}{5}{1} & = & \ciaa{2}{5}{} + \cia{1}{2}{5}    &
& \ciaa{1}{4}{} + \ciaa{1}{3}{4} & = & \ciaa{1}{3}{} + \ciaa{1}{4}{3}\\
  \ciaa{1}{2}{} + \cia{1}{5}{2} & = & \ciaa{1}{5}{} + \cia{1}{2}{5}    &
& \ciaa{3}{4}{} + \ciaa{1}{4}{3} & = & \ciaa{1}{4}{} + \ciaa{3}{4}{1}\\
  \ciaa{1}{5}{} + \cia{2}{5}{1} & = & \ciaa{2}{5}{} + \cia{1}{5}{2}&
& \ciaa{3}{4}{} + \ciaa{1}{3}{4} & = & \ciaa{1}{3}{} + \ciaa{3}{4}{1}\\
\\
  \ciaa{1}{3}{} + \cia{3}{5}{1} & = & \ciaa{3}{5}{} + \cia{1}{3}{5}    &
& \ciaa{1}{4}{} + \cia{1}{5}{4} & = & \ciaa{1}{5}{} + \cia{1}{4}{5}\\
  \ciaa{1}{5}{} + \cia{3}{5}{1} & = & \ciaa{3}{5}{} + \cia{1}{5}{3}    &
& \ciaa{4}{5}{} + \cia{1}{5}{4} & = & \ciaa{1}{5}{} + \cia{4}{5}{1} \\
  \ciaa{1}{5}{} + \cia{1}{3}{5} & = & \ciaa{1}{3}{} + \cia{1}{5}{3}&
& \ciaa{4}{5}{} + \cia{1}{4}{5} & = & \ciaa{1}{4}{} + \cia{4}{5}{1}\\
\\
  \ciaa{2}{4}{} + \ciaa{2}{3}{4} & = & \ciaa{2}{3}{} + \ciaa{2}{4}{3} &
& \ciaa{2}{3}{} + \ciaa{2}{5}{3} & = & \ciaa{2}{5}{} + \ciaa{2}{3}{5}\\
  \ciaa{3}{4}{} + \ciaa{2}{4}{3} & = & \ciaa{2}{4}{} + \ciaa{3}{4}{2} &
& \ciaa{2}{5}{} + \ciaa{3}{5}{2} & = & \ciaa{3}{5}{} + \ciaa{2}{5}{3}\\
  \ciaa{3}{4}{} + \ciaa{2}{3}{4} & = & \ciaa{2}{3}{} + \ciaa{3}{4}{2}&
& \ciaa{3}{5}{} + \ciaa{2}{3}{5} & = & \ciaa{2}{3}{} + \ciaa{3}{5}{2}\\
\\
  \ciaa{2}{4}{} + \ciaa{4}{5}{2} & = & \ciaa{4}{5}{} + \ciaa{2}{4}{5} &
& \ciaa{3}{4}{} + \ciaa{4}{5}{3} & = & \ciaa{4}{5}{} + \ciaa{3}{4}{5}\\
  \ciaa{2}{5}{} + \ciaa{4}{5}{2} & = & \ciaa{4}{5}{} + \ciaa{2}{5}{4} &
& \ciaa{3}{4}{} + \ciaa{3}{5}{4} & = & \ciaa{3}{5}{} + \ciaa{3}{4}{5}\\
  \ciaa{2}{5}{} + \ciaa{2}{4}{5} & = & \ciaa{2}{4}{} + \ciaa{2}{5}{4}&
& \ciaa{3}{5}{} + \ciaa{4}{5}{3} & = & \ciaa{4}{5}{} + \ciaa{3}{5}{4}\\
\end{array}
$$}

\end{document}